\documentclass[12pt,epsfig]{amsart}
\usepackage{tikz}
\usepackage{hyperref}
\newtheorem{theorem}{Theorem}[section]
\newtheorem{prop}{Proposition}[section]
\newtheorem{cor}[theorem]{Corollary}
\newtheorem{definition}[theorem]{Definition}
\newtheorem{notation}[theorem]{Notations}
\numberwithin{equation}{section}

\begin{document}
	\title[ Interpolation]{Interpolation by maximal  and minimal surfaces}
	\author{Rukmini Dey*}
	\address{*International Centre for Theoretical Sciences,
		Bengaluru- 560 089,  Karnataka, India.}
	\email{rukmini@icts.res.in}

	\author{Rahul Kumar Singh**}
	\address{**Department of Mathematics, Indian Institute of Technology Patna, Bihta, Patna-801106, Bihar, India}
	\email{rahulks@iitp.ac.in}
	
	
	\maketitle
	
	\begin{abstract}
		
		In this article, we use the inverse function theorem for Banach spaces to interpolate a given real analytic spacelike curve $a$ in Lorentz-Minkowski space $\mathbb{L}^3$ to another real analytic spacelike curve $c$, which is ``close" enough to $a$ in a certain sense by constructing a maximal surface containing them.  Next we apply the same method to interpolate two   given real analytic curve $a$ in Euclidean space $\mathbb{E}^3$ and a  real analytic curve $c$, which is also ``close" enough to ``a" in a certain sense with a minimal surface. 
		Throughout this study, the Bj\"orling problem and Schwarz's solution to it play pivotal roles.
	\end{abstract}
	
	\section{\bf{Introduction}}
	A maximal surface in ${\mathbb L}^3$  is a spacelike surface with vanishing mean curvature at every point. The study of maximal surfaces has been extensive in the past few decades, mainly focusing on methods to generate examples of such surfaces (see for instance  ~\cite{AL, Lo1, LTW} etc.). The Plateau's problem for maximal surfaces pertains to the existence of a maximal surface bounded by a simple closed spacelike curve. Similarly, we consider the Plateau's problem for two non-intersecting simple closed spacelike curves, separated by a certain distance, and seek a maximal surface with these two curves as its boundary. We refer to this problem as the interpolation problem using maximal surfaces.
	
	When the given spacelike curves consist of circles lying in parallel planes, separated by some distance,     Akamine and  L\'{o}pez obtain an elliptic catenoid whose boundary corresponds to the given circles (\cite{AL}). In other words, we can say that the elliptic catenoid interpolates the given circles. 
	It is important to note that solutions to these problems are already known for minimal surfaces in ${\mathbb  E}^3$  (see for instance ~\cite{N}).

	In general,  finding a maximal surface that can interpolate between given spacelike curves is a challenging task. The problem we have addressed here can be summarized as follows: Given two spacelike curves in $\mathbb{L}^3$ that are ``close" to each other in a certain sense (i.e., one curve belongs to a certain neighborhood of the other), we demonstrate the existence of a maximal surface that contains both curves.
	
	In a previous work (\cite{RPR}),  Dey,  Kumar and Singh  solved the problem of interpolating a given spacelike closed curve and a point with a maximal surface in $\mathbb{L}^3$, ensuring that the point becomes a singularity on the surface. Other instances of interpolating curves by maximal surfaces have been studied by Fujimori, Rossman, Umehara, Yang, and Yamada (\cite{FRUYY}).
	
	Solving the Plateau problem for minimal surfaces has been approached through various methods, as demonstrated by the works of Douglas and Rado, (\cite{douglas, douglas2, rado}) and many other mathematicians.  In our case, we are interested in utilizing the inverse function theorem for Banach spaces to tackle the interpolation problem for both maximal and minimal surfaces.
	
	I.N. Vekua in (\cite{vekua}) used the implicit function theorem for Banach spaces to solve Plateau's problem (for minimal surfaces) for curves lying sufficiently close to a given plane curve in $\mathbb{E}^3:=(\mathbb R^3, dx^2+dy^2+dz^2)$.
	
	Hamilton in (\cite{H}), specifically in example $5.4.13$, addressed Plateau's problem for minimal surfaces when the surface is represented as a graph. This was achieved using the Nash-Moser inverse function theorem.
	
	It is worth noting that Plateau's problem for maximal surfaces in pseudo-hyperbolic spaces, with boundary data at infinity, has been considered by Labourie, Toulisse, and Wolf in (\cite{LTW}).
	
	Additionally, L\'opez studied the existence of spacelike constant mean curvature surfaces spanning two circular contours (\cite{Lo1}).

	In this article, we utilize the inverse function theorem for Banach spaces to interpolate a given real analytic spacelike curve $a$ in $\mathbb{L}^3$ to another real analytic spacelike curve $c$ (which is sufficiently ``close" to $a$ in a certain sense) with  a maximal surface.  The maximal surface $X$ is defined as the real part of a specific complex analytic map from a complex domain to $\mathbb{C}^3$. The map $X$ is given by $X = \text{Re}(a+ id)$, where $a,d:\Omega\subset \mathbb{C} \rightarrow \mathbb{C}^3$ are complex analytic maps, and they satisfy the condition $\langle (a+id)^{\prime}, (a+id)^{\prime}\rangle_L =0$. This condition indicates that the map $a+id$ defines an isotropic curve in $\mathbb{C}^3$, and as a consequence, $X$ becomes maximal.
	
	Here, $\langle \cdot, \cdot\rangle_L$ represents the extended complex Lorentzian inner product (as described in the following paragraph). The solution to the Bj\"orling's problem is of significant importance in this context.
	
	\textit{Bj\"orling problem:}
	Recall that ${\mathbb L}^3 $ is defined as $({\mathbb R}^3, \langle \cdot, \cdot\rangle_L)$, where $\langle v, w\rangle_L = v_1 w_1 + v_2 w_2 - v_3 w_3$ for $v = (v_1, v_2, v_3) \in {\mathbb R}^3$ and $w=(w_1, w_2, w_3) \in {\mathbb R}^3$. The symbol $\times_L$ denotes the cross product of ${\mathbb L}^3$. Let $a: I\to \mathbb L^3$ be a real analytic spacelike curve as defined in Notation \ref{notation1}, with $a'(u)\neq 0$ for all $u\in I$.
	
	Also, consider $n: I\to \mathbb L^3$ as a real analytic map such that $ \langle n(u), a'(u) \rangle_L =0$ for all $u\in I $ and $n_1^2 (u)+ n_2^2(u) - n_3^2(u) = -1$ for all $u\in I$, where $n = (n_1, n_2, n_3)$. The map $n$ is referred to as the unit normal vector field along $a$, and it is timelike. The Bj\"orling problem poses the question of whether there exists a maximal surface $X : \Omega \to\mathbb L^3$ with $I\subset \Omega\subset \mathbb{C}$, where $\Omega$ is simply connected, and the surface $X$ satisfies the following conditions: $X(u, 0) = a(u)$ and $N(u, 0) = n(u)$ for all $u \in I$. Here, $N: \Omega\to \mathbb L^3$ is a normal vector field to the surface $X$. The solution to Bj\"orling's problem is given by (\cite{Mi})
	\begin{equation}
		X(u,v)=\text{Re}\{a(w)+i\int_{u_0}^{w} n(\widetilde{w})\times_L a'(\widetilde{w})d\widetilde{w}\};; w=u+iv, u_0\in I,
	\end{equation}
	where $n(w)$ and $a(w)$ are the complex analytic extensions of $n(u)$ and $a(u)$ over $\Omega$.

	\textit{Interpolation problem:} As mentioned earlier, the main objective of our study is to interpolate two closely related spacelike curves by leveraging the solution to the Bj\"orling problem. In essence, we aim to fix one curve and construct an appropriate neighborhood for it. Then, we seek to find another spacelike curve within this neighborhood and establish suitable Bj\"orling data along the first curve, (which includes a choice of a normal). With these conditions met, we can construct a Bj\"orling maximal surface that interpolates both curves. To state this more precisely, let $a: I\rightarrow \mathbb L^3$ and $l: I\rightarrow \mathbb L^3$ be two real analytic curves. We extend these curves analytically to a common domain $\Omega$, denoted as $a,l:\Omega\rightarrow \mathbb{C}^3$. By drawing inspiration from the solution to the Bj\"{o}rling problem, we pose the question: Can we find maps $d:\Omega\rightarrow \mathbb{C}^3 $ and $\gamma:\Omega \rightarrow \Omega$ (an automorphism) such that $Re((a+id)(u))=a(u)$ and $Re((a+id)(\gamma(u))=l(u)$? If we succeed, we have a generalized maximal surface of Bj\"{o}rling type given by $X:=Re(a+id)$, which serves as an interpolation between $a$ and $l$.
	
	We have succeeded in answering this question for  maximal and minimal surfaces. 
	
	The layout of the paper is as follows:
	
	In section \ref{IBS}, we invert the  Bj\"orling formula.  In section \ref{IMS} we fix a spacelike curve $a \in {\mathbb L}^3$ and
	construct a ``neighbourhood'' such that if another spacelike curve $l \in {\mathbb L}^3$ lies in that neighbourhood, there exists a maximal  surface which contains both $a$ and $l$.  As a corollary to the main theorem we show that the same method holds for minimal surfaces (the norms and cross product are  Euclidean norms and the Bj\"orling formula for minimal surface has to be used).  There are some subtle differences too as we elaborate later. Thus, in section \ref{IMiS} we get an interpolation by minimal surfaces in ${\mathbb E}^3$ as well.

	\section{\textbf{Inverting the Bj\"{o}rling-Schwarz  Function}}\label{IBS}
	
In this section we  will introduce and fix some notations that will be used throughout this article.

	The main aim in this section is to invert the Bj\"{o}rling-Schwarz  function $$ H:C^{\omega}(\Omega,\Omega)\times  C^{\omega}( \Omega, {\mathbb C}^3) \rightarrow   C^{\omega}(\Omega, {\mathbb C}^3)$$ defined by 
		\begin{equation}\label{defnofH}
			H(\gamma, d)= (a+id)\circ \gamma
		\end{equation} (see Defintion \ref{defnofH}).  In the process we  define certain space ${\mathcal J} \subset C^{\omega}(\Omega, {\mathbb C}^3)$  which is the  union of slices $J_n$ (see Defintion \ref{notation2}) and invert the map $H $ (see Proposition \ref{defnofH}).   For this we compute the differential $DH$ of $H$,  restrict to $T_n$, i.e.,  $DH|_{T_n}$ (see Definition \ref{notation3}).  We define ${\mathcal H}_n$ (see Defintion \ref{notation3}) and show later that it is the same as    $DH|_{T_n}$.  We show that ${\mathcal H}_n: T_n  \rightarrow J_n$ is invertible.  Then invoking the inverse function theorem for Banach spaces we show that the 
Bj\"{o}rling-Schwartz  function $H$ is invertible in suitable domains.

	\begin{notation}\label{notation1} We introduce the following notations:
	
		Let $\Omega$ denote the open unit disk in the complex plane, and for simplicity, let $\gamma_0(z) = z$ where $\gamma_0:\Omega \rightarrow \Omega$.
		\begin{enumerate}
			\item We denote the Euclidean and Lorentzian inner products as $\langle\cdot , \cdot \rangle_E$ and $\langle \cdot, \cdot\rangle_L$, respectively. Moreover, we denote the perpendicular operation with respect to $\langle \cdot , \cdot \rangle_E$ as $\perp_E$, and with respect to $\langle \cdot , \cdot \rangle_L$ as $\perp_L$. Additionally, we use $\times_E$ to represent the cross-product in the Euclidean sense, and $\times_L$ for the Lorentzian sense.
			\item Let $a: I\to\mathbb R^3$ be a real analytic curve that can be analytically extended to $\Omega$ containing $I$. This extension is also differentiable on $\overline{\Omega}$. Furthermore, for all $z\in \overline{\Omega}$, we assume $a^\prime(z)\neq 0$, where $\prime$ as usual denotes the derivative.
			
			\item For a bounded function $f:\overline{\Omega}\rightarrow \mathbb{C}^n$, where $f=(f_j)$ and $j=1,2,...,n$, we define its norm as $\|f\|_{\infty}:=\max\limits_{j}\{\sup_{z\in \overline{\Omega}}|f_j(z)|\}$, also known as the sup norm.
			
			\item For $n\geq 1$, let $C^{\omega}(\Omega, \mathbb{C}^n)$ denote the space of all complex analytic maps on $\Omega$, continuously extendable to the boundary of $\Omega$ (i.e., $\overline{\Omega}$). This set forms a Banach space under the norm $\|f\|_{\infty}$.
			
			\item Let $C^{\omega}(\Omega, \Omega)$ denote the space of $\Omega$-valued complex analytic functions defined on $\overline{\Omega}$, which also function as automorphisms of $\Omega$. This space is an open subset of $C^{\omega}(\Omega, \mathbb{C})$.
			
			\item For a fixed $u_0\in I$, we define the following:
			\begin{enumerate}
				\item $\mathcal C_{(\Omega,\Omega, u_0)}:=\{ \gamma \in C^{\omega} ( \Omega, \Omega) \mid \gamma(u_0) = u_0 \}$
				\item ${\mathcal C}_{(\Omega,\mathbb C,u_0)}:=\{ V  \in C^{\omega}(\Omega, {\mathbb C} ) \mid V(u_0) =0\}$ 
				\item $\mathcal W_{(\Omega,\mathbb C^3)}:=\{s  \in C^{\omega}(\Omega,\mathbb{C}^3) \mid s(u) \text{ is real on } u \in I\}$. This also forms a Banach space.
				\item $\mathcal W_{(\Omega,\mathbb C^3, u_0)}:=\{s  \in   C^{\omega}(\Omega,\mathbb{C}^3) \mid s(u) \text{ is real for all }  u \in I\; \text{ and } s(u_0) = 0\}$. This is a closed subset of $C^{\omega}(\Omega,\mathbb{C}^3) $, and hence is a Banach space with the induced norm.
				\item ${\mathcal L_{(I,\mathbb L^3)}}:= \{c\in C^{\omega}(\Omega,\mathbb{C}^3)  \mid c \;\text{and}\; c' \;\text{are real on}\; I\}$. 
				\item ${\mathcal S}_p:=\{c\in \mathcal L_{(I,\mathbb L^3)} \mid  \langle c^\prime, c^\prime \rangle_{L} >0 \text{ on } I\}$.  
				\item $\widetilde{\mathcal  S_p}:=\{ c \in C^{\omega}(\Omega,\mathbb{C}^3) \mid \text{Re}(c^{\prime}(u) ) \text{ is spacelike, non-zero} \} $.
				
				
				\item ${\mathcal N}:=\{n\in C^\omega(\Omega, \mathbb C^3) \mid n|_I \text{ is real, } \langle n, n\rangle_L=-1,\; \langle n,  a^{\prime}\rangle_L =0\}$. This set is closed in the subspace sup-norm topology. 
			\end{enumerate}
		\end{enumerate}
	\end{notation}

	\begin{prop} We have the following proposition. 	
		\begin{enumerate}
			\item[(a)] ${\mathcal S}_p$  is an open subset of ${\mathcal L_{(I,\mathbb L^3)}}$ in subspace  sup-norm topology. 
			\item[(b)] $\widetilde{{\mathcal  Sp}}$ is an open subset of $C^{\omega}(\Omega,\mathbb{C}^3)$ in subspace sup-norm topology. 
		\end{enumerate}
	\end{prop}
	\begin{proof}
		(a). 
		We take the map $g : \Omega \times {\mathcal L_{(I,\mathbb L^3)}} \rightarrow {\mathbb C}$  given by $g(z , c) = c_{1}^{\prime 2}(z) + c_{2}^{\prime 2}(z) - c_{3}^{\prime 2}(z)$.  It is a  continuous map. Therefore,  by restriction,  the map $\phi = g : I \times {\mathcal L_{(I,\mathbb L^3)}} \rightarrow {\mathbb R}$ is continuous.  Then $\phi^{-1} ({\mathbb R}^+) =  I \times {\mathcal S}_p$ is an open set in $I \times {\mathcal L_{(I,\mathbb L^3)}} $ since  $\phi$ is continuous and ${\mathbb R}^+$ is an open subset of ${\mathbb R}$.  It follows easily that ${\mathcal S}_p$ is an open subset of ${\mathcal L_{(I,\mathbb L^3)}}$.   
		
		(b).  Similar to case (a), we take the continuous map $f: \Omega \times  C^{\omega}(\Omega,\mathbb{C}^3  ) \rightarrow {\mathbb C}$  given by $f(z, c) = || \rm{Re}(c^{\prime}(z))||_L^2 $.  We restrict  $f$ to $I \times C^{\omega}(\Omega,\mathbb{C}^3)$ and name this map $\phi$.  Thus $\phi$ maps  $I \times C^{\omega}(\Omega,\mathbb{C}^3)$ to ${\mathbb R}$ and is continuous.  Then  $\phi^{-1} ({\mathbb R}^+) =  I \times \widetilde{\mathcal Sp}$ is an open set in $I \times  C^{\omega}(\Omega,\mathbb{C}^3)$ since  $\phi$ is continuous and ${\mathbb R}^+$ is an open subset of ${\mathbb R}$.  It follows that $\widetilde{{\mathcal  Sp}}$ is an open subset of $C^{\omega}(\Omega,\mathbb{C}^3)$. 
	\end{proof}

	{\bf Note:} 
	We use the notion of smooth maps and derivations of maps  between Banach spaces in the same sense as in  Kriegel and Michor  (\cite{KM}).

	\begin{prop}\label{differential}
		
		Given $a \in {\mathcal S_p} $, the map $ H:C^{\omega}(\Omega,\Omega)\times  C^{\omega}( \Omega, {\mathbb C}^3) \rightarrow   C^{\omega}(\Omega, {\mathbb C}^3)$ defined by 
		\begin{equation}\label{defnofH}
			H(\gamma, d)= (a+id)\circ \gamma
		\end{equation} is smooth.  
		Let  $  \gamma_0\in C^{\omega}(\Omega,\Omega)$, $ \gamma_0(z)= z$  and $ d^0 \in  C^{\omega}(\Omega, {\mathbb C}^3)$. Then the derivative  of $H$ at the point $(\gamma_0,d^0)$, $ DH_{(\gamma_0,d^0)}:C^{\omega}(\Omega, \mathbb{C})\times  C^{\omega}(\Omega, {\mathbb C}^3)  \rightarrow  C^{\omega}(\Omega, {\mathbb C}^3)     $ is given by 
		
		$DH_{(\gamma_0,d^0)}(V, d) = ((a^{\prime} + i {d^0}^{\prime})V  +i d)$.
	\end{prop}

	\begin{proof}

		One can  show that the map $ H(\gamma,d)=  (a+id)\circ \gamma$ is smooth.

		Let $ P_0=(\gamma_0, d^0) $ and $ h=(V, d)$, where $V\in C^{\omega}(\Omega,\mathbb C)$ and $d\in  C^{\omega}(\Omega, {\mathbb C}^3)$. Then
		\begin{eqnarray*} DH_{P_0}(h)&=& \lim_{s\rightarrow 0}\frac{1}{s}\{H(P_0+sh)-H(P_0)\}\\
			&=& \lim_{s\rightarrow 0}\frac{1}{s}\{H(\gamma_0+sV, d^0+sd)-H(\gamma_0,d^0)\}\\
			&=& \lim_{s\rightarrow 0}\frac{1}{s}\{(a+i(d^0+sd))\circ(\gamma_0+sV)-(a+id^0)\circ\gamma_0\}.\\
		\end{eqnarray*}
		Let $a=(a_1,a_2,a_3)$ and $d=(d_1^0,d_2^0,d_3^0)$. Then we have 
		$$\lim_{s\rightarrow 0} \frac{1}{s}\{a_j(\gamma_0(z)+sV(z))-a_j(\gamma_0(z))\}
		= V(z)a_j'(\gamma_0(z))\; ; j=1,2,3$$ and
		\begin{eqnarray*}
			\lim_{s\rightarrow 0}\frac{1}{s}\{(d_j^0+sd_j)(\gamma_0+sV)-d_j^0(\gamma_0)\}
			&=& V(z)d_j^0{'}(\gamma_0(z))+d_j(\gamma_0(z))~~\text{where}~~j=1,2,3.\\
		\end{eqnarray*}
		Therefore, we must have
		\begin{equation*}
			DH_{(\gamma_0,d^0)}(V,d)=\{z\longmapsto a'(z)V(z)+iV(z)d^{0'}(z)+id(z)\}.
		\end{equation*}
		Or in other words, we have	
		\begin{equation*}
			DH_{(\gamma_0,d^0)}(V,d)=   (a'+i{d^0}^{\prime})V+ i d.
		\end{equation*}
	\end{proof}

	Let $a$ be as before and $n \in\mathcal{N}$.  On $\Omega$, we define 
	\begin{equation}\label{d0n}
		d_0^n(z)= \int_{u_0}^z  (n \times_L a^{\prime})(w) dw.
	\end{equation}
	We see $d_0^n(u_0) =0$, and its restriction on $I$ is space-like. Therefore  $d_0^n\in \mathcal{W}_{(\Omega, \mathbb C^3, u_0)}$.  
	
	Moreover, $$\langle d_0^{n \prime}, d_0^{n \prime} \rangle_L = \langle a^{\prime}, a^{\prime}\rangle_L,\;\; \langle d_0^{n \prime}, a^{\prime}\rangle_L =0.$$    
	In other words, $a+i d_0^n$ is L-isotropic in the sense that $\langle (a+i d_0^n)^{\prime}, (a+id_0^n)^{\prime}\rangle_L =0$.	We define the function $\Delta^n: I\to \mathbb R^3$ as 
	\begin{equation}\label{deltan}
		\Delta^n(u) := d_0^{n \prime}(u) \times_E a^{\prime}(u).
	\end{equation}
	{\bf Note:} Since $d_0^{n \prime}(u)$ and $a^{\prime}(u)$ are spacelike,  $\Delta^n(u)$ is timelike.
	\begin{definition}\label{notation2}  We have the following
		\begin{enumerate}
			\item $J_n:= \{ c \in  C^{\omega}(\Omega, {\mathbb C}^3)| \text{Re} (c^{\prime}(u))  \perp_E \Delta^n(u)  \; \forall \; u \in I  \; \rm{and} \; \text{Re} (c^{\prime}(u)) \; \text{is} \; \text{space-like} \}$.  
			
			\item ${\mathcal J} := \cup_{n  \in \mathcal N} J_n$ with the topology as follows:  $V\subset C^{\omega}(\Omega, {\mathbb C}^3)$ is open in ${\mathcal J}$ if and only if for all $n$, $V \cap J_n$ is open in $J_n$ in sup-norm topology. We  call this $n$-topology.  
		\end{enumerate}

	\end{definition}
	Since $\mathcal{J}\subset C^\omega(\Omega, \mathbb C^3)$, it has a sup norm topology as given in Notations \ref{notation1}. By notation $\mathcal{J}$, we always mean $n$-topology unless otherwise stated. We discussed the relation of this topology with the induced sup norm topology  on $\mathcal{J}$ later in the Proposition ~\ref{nbhd}.
	
	Recall $H(\gamma, d)=(a+id)\circ \gamma$, $H(\gamma_{0}, 0)=a$. It is easy to see that $a\in J_n$ for all $n$. Therefore, $H^{-1}(J_{n})$ is non-empty for all $n$.

	\begin{definition}\label{notation3} We define the following-
		\begin{enumerate}
			\item For $H$ as in the equation \eqref{defnofH} and $J_n$ as in the Definition \ref{notation2}, we define $X_n:= H^{-1}(J_n)$, and  $H_n:= H|_{X_n}$.  The co-domain of the  map $H_n$ is  $J_n$. We have  $H_n(\gamma, d) = (a+id) \circ \gamma$.
			\item For each $n\in \mathcal{N}$, we define  
			$ T_{n}:= \{ (V, d) \in {\mathcal C}_{(\Omega,\mathbb C, u_0)} \times {\mathcal W_{(\Omega,\mathbb C^3, u_0)}} $: $\text{Re} ((a^{\prime}(u) + i d^{n \prime}_0(u)  ) V(u) + i d(u))^{\prime} \perp_E \Delta^n(u),\;u\in I$\}. 
			
			\item Let ${\mathcal H}_n : T_n \rightarrow J_n  $ be a map defined by ${\mathcal H}_n(V, d) = (a^{\prime} + i d^{n \prime}_0 ) V + i d$. The $ Re((a^{\prime} + i d^{n \prime}_0) V + i d)$ is linear combination of $a'$ and $d_0^{n\prime}$ on $I$. Thus, ${\mathcal H}_{n}(T_n)\subset J_n$.
		\end{enumerate}
	\end{definition}
	
	In the below, we aim to relate $H_n$ and $\mathcal{H}_n$, and as one can guess, $\mathcal{H}_n$ turns out to be the differential of $H_n$ at some point.  We will prove it. 
	
	We start analyzing the domain $H_n^{-1}(J_n)= H^{-1}(J_n)$. 
	\begin{prop}\label{Bmnfld}
		There is an open set ${\mathcal U}_n \subset H_n^{-1}(J_n)$ containing $p_0= (\gamma_0, d_0^n)$ which  has a Banach manifold structure.
	\end{prop}
	\begin{proof}
		
		Let $E_{\epsilon} = \{ (V,d) \in {\mathcal H}_n^{-1}(J_n) : \gamma_0 + \epsilon V \; \text{is} \; \text{invertible} \}$. 
		\begin{enumerate}
			\item For some $\epsilon$, $E_\epsilon$ may be empty. 
			\item if for some $\epsilon$, $E_\epsilon\neq \phi$, then it is an open subset of $T_n$.  
			
			\item  if for some $\epsilon_0$, $E_{\epsilon_0}\neq \phi$, then for every $\epsilon<\epsilon_0$, $E_{\epsilon}\subset E_{\epsilon_0}$, and $E_{\epsilon}\neq \phi$. 
		\end{enumerate}
		
		Choose $\epsilon_0$ such that $E_{\epsilon_0}\neq \phi$,  $$E_n = \cup_{0<\epsilon<\epsilon_0} E_{\epsilon}$$ is open in ${\mathcal H}_n^{-1}(J_n)$. 
		
		For $0<\epsilon < \epsilon_0$ we take an open subset of $H_n^{-1}(J_n )$ as follows
		$${\mathcal \theta}_{\epsilon}:=\{ (\gamma_1,  d_1) \in H_n^{-1}(J_n)|:\;  \;  \|\gamma_1-\gamma_0\|_\infty < \epsilon,  \|d_1-d_0^n\|_\infty < \epsilon \}.$$ Let ${\mathcal U}_n = \cup_{0<\epsilon<\epsilon_0} \theta_{\epsilon}$. This is open in the subspace topology of $H_n^{-1}(J_n)$.
		
		We define  $\psi_{\epsilon} : \theta_{\epsilon} \rightarrow E_{\epsilon} $ by the map 
		$$\psi(\gamma_1, d_1)=\left(\frac{\gamma_1-\gamma_0}{\epsilon},  \frac{d_1-d_0^n}{\epsilon}\right)= (V, d).$$ 
		
		To see that the image $(V, d)$ is indeed in $T_n$, we show that $\text{Re} ((a+id_0^n)^{\prime}  V + i d)^{\prime} )\perp_E \Delta^n$ on $I$.

		Let $c = (a + i (d_0^n + \epsilon d)) \circ (\gamma_0 + \epsilon V)$.  Since for all $0 < \epsilon < \epsilon_0$, $(\gamma_1, d_1) = ( \gamma_0 + \epsilon V, d_0^n + \epsilon d)$ belongs to $H^{-1}(J_n)$,  $\text{Re}(c^{\prime}(u)) \perp_E \Delta^n$ on $I$. 
		
		On the other hand we see that on $I$, $\text{Re} (c^{\prime}) =  \text{Re} (a + i (d_0^n + \epsilon d) )^{\prime} \circ (\gamma_0 + \epsilon V) \cdot (1 + \epsilon V^{\prime}))$.  Therefore for all $\epsilon<\epsilon_0$, $\text{Re} ((a + i (d_0^n + \epsilon d) )^{\prime} \circ (\gamma_0 + \epsilon V) \cdot (1 + \epsilon V^{\prime}))$ is Euclidean perpendicular to $\Delta^n$.
		
		Moreover $\text{Re} \frac{d c^{\prime}(u)}{d \epsilon}|_{\epsilon \rightarrow 0} =\text{Re} (\frac{d c(u)}{d \epsilon}|_{\epsilon \rightarrow 0})^{\prime}$, and $\frac{d c(u)}{d \epsilon}|_{\epsilon \rightarrow 0} = (a + i d_0^n)^{\prime}(u) V(u) + i d(u) $. Therefore  $\text{Re} ((a + i d_0^n)^{\prime}(u) V(u) + i d(u))^{\prime} \perp \Delta^n $ on $I$.

		The map $\psi$ extends to a map $\psi: {\mathcal U}_n \rightarrow E_n$ since $\theta_{\epsilon}$ and $E_{\epsilon}$ are nested sets respectively. An inverse of $\psi$ also exists. 
		
		We take $(\mathcal{U}_n, \psi)$ as a Banach manifold modelled on $E_n$. 
	\end{proof}
	\color{black}
	For $p_0 = (\gamma_0, d_0^n)$, we claim the tangent space of $\mathcal{U}_n$ at $p_0$ is  $T_n$ as given in the Definition \ref{notation3}. 
	
	\begin{prop}\label{tangent}
		For $\mathcal{U}_n$ as in the Proposition \ref{Bmnfld},  and $T_n$ as in the Definition \ref{notation3},  we have  $$T_{p_0} (\mathcal{U}) = T_n.$$
	\end{prop}	
	
	If we elaborate above equality, we are claiming that $T_{p_0}(\mathcal{U})$  consists of all $(V, {D}) \in \mathcal C_{(\Omega,\Omega, u_0)}\times {\mathcal W}_{(\Omega,\mathbb C^3, u_0)} $ such that on $I$, $\text{Re} ((a^{\prime}(u) + i d^{n \prime}_0(u)  ) V(u) + i {D}(u))^{\prime} \perp_E \Delta^n(u) $. 
	\begin{proof}
		Take a curve $c: (-\epsilon, \epsilon)\to \mathcal{U}\subset H_n^{-1}(J_n)$ such that $c(0)= (\gamma_0, d_0^n)$.   By Taylor series expansion,   $c(t)  = (\gamma(t), d(t))=(\gamma_0 + t V + o(t^2),  d^n_0 + t {D} + o(t^2))$.   We denote $D_t= H_n(c(t))\in J_n$,  we have
		\begin{equation}\label{l1}\text{Re}(D^{\prime}_t(u))\perp_E \Delta^n(u) \text{ on } I.\end{equation}   Note that $D_t^\prime(u)= \frac{d}{du} D_t(u)$. Moreover, a similar calculation as in Proposition  \ref{differential} shows that $\frac{d D_t}{dt}|_{t=0}  =  (a^{\prime} + i d^{n\prime}_0  ) V + i {D}$.

		Differentiating the above equation \eqref{l1} at  $t=0$, we get 
		$\langle Re(\frac{d}{dt}D_t^\prime(u)|_{t=0}), \Delta^n\rangle_{E}=0$, and we know $\text{Re}(\frac{dD^{\prime}_t}{dt}|_{t=0}) =\text{Re} (\frac{d D_t}{dt}|_{t=0})^{\prime}$. Therefore $\text{Re} (\frac{d D_t}{dt}|_{t=0})^{\prime} \perp_E \Delta^n$. 
		
		This proves 
		$$\text{Re} ((a^{\prime} + i d^{n\prime}_0  ) V + i {D})^{\prime}) \perp_E \Delta^n;\; \text{that is } T_{p_0}(\mathcal{U})\subset T_n.$$
		
		Other way is direct, as for given $(V,D)\in T_n$, we construct $c(t)= (\gamma_0+tV, d_0^n+tD)$.  For small $t$, such $c$ is a curve in $\mathcal{U}$.

	\end{proof}

	For $p_0$, $\mathcal{U}$ as before, we have the following proposition, and $H_n: \mathcal{U}\subset X_n\to J_n$,  we have 
	$D(H_n)_{(\gamma_0, d_0^n)}: T_{p_0}\mathcal{U}= T_n\to J_n$ turns out be $\mathcal{H}_n$ as in the Definition \ref{notation3}.

	Below, we will show that $D(H_n)_{(\gamma_0, d_0^n)}(= {\mathcal H}_n)$ is invertible.   
	
	\begin{prop}\label{invertible}
		The map $ {\mathcal H}_n : T_n   \rightarrow  J_n $ is invertible. 
	\end{prop}
	
	As many notations are already used, we recall a few here:   $\gamma_0 \in C^\omega(\Omega, \Omega)$  is defined as  $\gamma_0 (z) = z$ and $d_0^n \in {\mathcal W}_{(\Omega,\mathbb C^3, u_0)}$ as in the equation \eqref{d0n}. Here $p_0 = (\gamma_0, d_0^n) $.   Moreover $
	T_n$ is defined in (2) of Definition \ref{notation3} and $T_n \subset  \mathcal C_{(\Omega,\mathbb C, u_0)} \times     {\mathcal W}_{(\Omega,\mathbb C^3, u_0)}$
	\begin{proof}
		We know that
		${\mathcal H}_n(V, d)= (a'+id_0^{n'})V +id $, and since $(V,d)\in T_n$, $(V,d)$ is related by  $  \text{Re} ((a'+id_0^{n'})V +i{d}))^{\prime} \perp_E \Delta^n$ on $I$.

		We want to show  that ${\mathcal H}_n$ is invertible, i.e., given a fixed   $s =(s_1, s_2, s_3)\in J_n$ there exists a unique  $ (V, {d})  \in \mathcal C_{(\Omega,\Omega, u_0)} \times {\mathcal W}_{(\Omega,\mathbb C^3, u_0)}$ such that ${\mathcal H}_n (V, {d})(z) = s(z)$, i.e.  we have to solve for $(V,d)$ satisfying 
		$$ s(z) = (a'(z)+id_0^{n'}(z))V(z)+i{d}(z), \text{and} $$
		$$ \text{Re} ((a'+id_0^{n'})V +i{d}))^{\prime} \perp_E \Delta^n.  $$
		
		Since $s \in J_n$, $\text{Re} (s^{\prime}) \perp_E \Delta^n$ on $I$. Let $\tau(z) = s^{\prime}(z)$ and thus $\text{Re}(\tau(u)) \perp_E  \Delta_n$ on $I$. 
		
		As a first step, we want to solve for $\widetilde{V}$ and ${\widetilde{d}}$ such that   $\widetilde{V} \in C^\omega(\Omega, \mathbb{C}) $ and ${\widetilde{d}} \in \mathcal W_{(\Omega,\mathbb C^3)}$ and satisfying the 
		
		\begin{equation}\label{l3}
			\tau(z) = s^{\prime} (z) = (a'(z)+id_0^{n'}(z) ) \widetilde{V}(z) +i {\widetilde{d}}(z)).
		\end{equation}

		Let us denote $\tau= (\tau_1,\tau_2, \tau_3)$ and  $$ (A(z),B(z),C(z))=(a_1'+id_{01}^{n'} ,a_2'+id_{02}^{n'},a_3'+id_{03}^{n'})(z), $$ \
		$$ (X(z),Y(z),Z(z))=(i\{\widetilde{d}_1,i{\widetilde{d}}_2,i{\widetilde{d}}_3)(z), $$ and 
		$$\widetilde{V} = \widetilde{V}_1 + i \widetilde{V}_2.$$
		
		The equation \eqref{l3} gives set of three equations
		\begin{equation}\label{xyz2}
			\begin{split}
				(A(z)\widetilde{V}(z) + X(z))=\tau_1(z),\\
				(B(z)\widetilde{V}(z) + Y(z))=\tau_2(z),\\
				(C(z)\widetilde{V}(z) + Z(z))=\tau_3(z),
			\end{split}
		\end{equation}
		
		that is we have 	
		\begin{equation}\label{xyz3}
			\begin{split}
				X(z)=\tau_1(z) - A(z)\widetilde{V}(z) ,\\
				Y(z)=\tau_2(z) - B(z)\widetilde{V}(z) ,\\
				Z(z)=\tau_3(z)- C(z) \widetilde{V}(z).
			\end{split}
		\end{equation}
		
		We have an  additional constraint that ${\widetilde{d}}={\widetilde{d}}(u,0)$ is real on $I$.  
		Using this fact and from the the equations \eqref{xyz3}, we need $X, Y, Z$ and $\tilde{V}$ such that on  $u \in I$ we have 
		\begin{equation}\label{xyz4}
			\begin{split}
				0= -\text{Im} {\widetilde{d}}_1(u) = \text{Re} X(u)= \text{Re} (\tau_1(u) - A(u)\widetilde{V}(u)) ,\\
				0= -\text{Im} {\widetilde{d}}_2(u) = \text{Re} Y(u)=\text{Re}(\tau_2(u) - B(u) \widetilde{V}(u)) ,\\
				0= -\text{Im} {\widetilde{d}}_3(u) = \text{Re} Z(u)=\text{Re} (\tau_3(u)- C(u) \widetilde{V}(u)).
			\end{split}
		\end{equation}
		
		The vector $\Delta^{n}  = d_0^{n \prime} \times_E a^{\prime}   = (\Delta^n_{1}, \Delta^n_{2}, \Delta^n_{3})  $ is not zero on $ I$.  At $u \in I$, suppose without loss of generality, $ \Delta^n_{3}(u) =( a_2'd_{01}^{n \prime} -a_1'd_{02}^{n \prime})(u) \neq 0$. 
		
		At that $u \in I$, 
		writing $\widetilde{V}(u) = \widetilde{V}_1(u) + i \widetilde{V}_2(u)$ with $V_1$ and $V_2$ real, we have from the first two equations:

		\begin{equation}
			\begin{split}
				\widetilde{V}_1(u)=\frac{ \text{Re} \tau_2 d_{01}^{n\prime}- \text{Re} \tau_1 d_{02}^{n \prime}}{\Delta^n_3}(u),\\
				\widetilde{V}_2(u)=\frac{\text{Re} \tau_1 a_2'- \text{Re} \tau_2 a_1^{\prime}}{-\Delta^n_3}(u).
			\end{split}
		\end{equation} 
		
		There are two other ways of writing $\widetilde{V}_1$ and $\widetilde{V}_2$ on $I$ when $\Delta^n_1 \neq 0$ or $\Delta^n_2 \neq 0$.  But in at $u$ when all $\Delta_i^n(u)\neq 0$, all expressions turn out to be the same (see appendix \ref{Vconsistent}). 	
		
		This gives a well-defined function on $\tilde{V}$ on $I$. 
		
		Moreover,  such $\tilde{V}$ can be analytically extended on $\Omega$ (see appendix \ref{defn of tilde V on all of Omega}).

		Substituting $\tilde{V}$, we will get $X, Y,$ and $Z$.   Moreover we can verify that for such $\tilde{d}$,  ${\widetilde{d}}(u,0)$ is real on $I$ (see appendix \ref{d(u,0)realonI}). 
		
		Thus we have 
		\begin{equation}\label{es}
			\tau(z) = s^{\prime} (z) = (a'(z)+id_0^{n\prime}(z) ) \widetilde{V}(z) + i {\widetilde{d}}(z)
		\end{equation}
		with $\widetilde{V} \in C^\omega(\Omega, \mathbb{C}) $ and ${\widetilde{d}} \in \mathcal W_{(\Omega,\mathbb C^3)}$. In other words,  we have $s^{\prime} = (a^{\prime} + i d_0^{n \prime}) \widetilde{V} + i {\widetilde{d}}$.

		To complete the proof,  we need to solve for $(V, d) \in T_n \subset {\mathcal C}_{(\Omega,\Omega, u_0)} \times {\mathcal W}_{(\Omega,\mathbb C^3, u_0)}$ such that $s = (a^{\prime} + i d_0^{n \prime}) V + i d$.  That is we need to solve for $(V, {d})$ such that  	
		\begin{equation}\label{ODE}
			s^{\prime} = (a^{\prime} + i d_0^{n \prime}) \widetilde{V} + i {\widetilde{d}} = ((a^{\prime} + i d_0^{n \prime}) V + i {d})^{\prime}.
		\end{equation}
		
		We claim that there is a unique solution $ (V, {d}) \in T_n$  to the above equation. First, we take 
		${d(w)} = \int_{u_0}^w {\widetilde{d}}(\widetilde{w}) d \widetilde{w}$. 
		
		Then the equation for $V$ is $(a^{\prime} + i d_0^{n \prime}) \widetilde{V} = ((a^{\prime} + i d_0^{n \prime}) V)^{\prime}$, with $V(u_0) =0$.  This is an ODE that can be solved and has a unique solution (see appendix \ref{solving ODE}).   Moreover, such $(V,d)\in T_n$.
		
		To prove the uniqueness of the solution of equation ~\eqref{ODE}, we take two solutions $(V^1, {d}^1) $ and $(V^2, {d}^2) $ and let $\chi = V^1 - V^2$ and $d = {d}^1 - {d}^2$. Then $(\chi, d)$ satisfy the equation
		\begin{equation}
			((a^{\prime} + i d_0^{n \prime}) \chi + i d)^{\prime} =0
		\end{equation}
		with $\chi(u_0) =0$ and $d(u_0) =0$.
		
		Then $(a^{\prime} + i d_0^{n \prime}) \chi + i d = c_0$, where the constant $c_0=0$. 
		On $I$,   $(a^{\prime} + i d_0^{n \prime})(u) \chi(u) + i d(u) =0.$
		
		Equating real and imaginary parts to zero,  we get $a^{\prime}(u) \chi(u) =0$ and $d_0^{n \prime}(u) \chi(u) = -d(u)$.  Since at least one of the components of $a^{\prime}(u)$ is nonzero, $\chi(u) =0$ and hence $d(u) =0$. It follows that  on $I$, $\chi =0$ and $d =0$ on $I$.  So by Identity theorem $d=0$ on $\Omega$.  Thus uniqueness is proved.

	\end{proof}

	\begin{prop}\label{invofHN}
		There is an open neighborhood $U_n\subset \mathcal{U}$ of $p_0$, and open set $V_n\subset J_n$ containing $H_n(p_0)=x_0$, such that $H_n: U_n\to V_n$ is diffeomorphism. 
	\end{prop}
	\begin{proof}
		We know that
		$DH_n|_{p_0}(V,{d})=  (a'+id_0^{n'})V +i{d} $.  We also have $T_n = T_{p_0}(\mathcal{U})$ and ${\mathcal H}_n = DH_{n| p_0}$. Then $DH_{n |p_0}$ is invertible since ${\mathcal H}_n$ is invertible, by Proposition \ref{invertible}.  Therefore the  map $DH_{ n|p_0}:   T_{p_0}(\mathcal{U}) \rightarrow  J_n $ is invertible. 
		
		Therefore using the inverse function theorem of the Banach manifolds, we get the required open sets. 
	\end{proof}	
	
	Now the above process can be repeated for each $n\in \mathcal{N}$, and we define
	$U^\Omega:= \cup_n U_n$ and $V^\Omega:= \cup_n V_n$.   We see that $U^\Omega\subset C^\omega(\Omega,\Omega)\times C^\omega(\Omega,\mathbb C^3)$ and $V^\Omega\subset \mathcal{J}\subset C^\omega(\Omega, \mathbb C^3)$.  We recall here that the topology of $\mathcal{J}$ is n- topology. 
	Now we will write the following corollary. 
	\begin{cor}\label{IFT}
		Given any $q\in V^\Omega$, there is a $p\in U^\Omega$ such that  $H(p)= q$.  In other words, given a $q\in V^\Omega$, there is a $\gamma\in C^\omega(\Omega,\Omega)$, and $d\in \mathcal{W}_{(\Omega, \mathbb C^3, u_0)}$ such that $(a+id)\circ\gamma= q$. 
	\end{cor}

	\begin{proof}
		Since $q\in V^\Omega$, there is a $n\in \mathcal{N}$, such that  $q\in V_n$.  By the proposition \ref{invofHN}, we have $(\gamma, d_n)\in U_n$, such that $H_n(\gamma, d_n)= q$. 
	\end{proof}

	Now we are ready to discuss the interpolation problem as defined in the introduction. 
	
	\section{\textbf{Interpolation by maximal surface}} \label{IMS}
	Let $a, d_0 \in \mathcal{S}_p \subset \mathcal{L}_{(I,\mathbb{L}^3)}$ and let  $x_0= a+ i d_0$ be such that $\langle a^{\prime}, a^{\prime}\rangle_L = \langle d^{\prime}_0,d^{\prime}_0\rangle_L$, and $  \langle a^{\prime}, d^{\prime}_0\rangle_L = 0$. We have $x_0 \in \mathcal{J}$.  Let $V^{\Omega}$ be an open subset of $\mathcal{J}$ in n-topology containing $x_0$.
	
	First,  in subsection $3.1$,  we will show that there is an open ball $B_\epsilon(x_0)$ in the sup-norm topology that lies in $V^{\Omega}$. This will help us to see the closeness of the two curves easily. Then in the subsection $3.2$ we prove the interpolation problem as discussed in the introduction.  
	

	\subsection{Existence of sup-norm open ball containing an appropriate $x_0$ in $V^{\Omega}$} 
	
	Before finding a sup-norm ball containing $x_0$, we prove that $\tilde{\mathcal{S}}_p \subset \mathcal{J}$.
	\begin{prop}\label{SptildemathcalJ}
		$\tilde{S}_p\subset \mathcal{J}$.
	\end{prop}
	\begin{proof}
		
		Let $s \in\tilde{\mathcal{S}p} $. We claim that $s \in J_{n_s}$ for some $n_s \in {\mathcal N}$. For this we just need to show for $u\in I,\;\text{Re}(s^{\prime}(u) ) \perp_E \Delta^{n_s}(u)$. 
		
		If $\text{Re}(s'(u))=\alpha(u) a'(u)$ for all $u \in I $, then $s\in J_n$ for all $n$, and we are done.
		
		Now, if the above is not true, then there will be some $u_0 \in I$, $\text{Re}(s'(u_0))\neq \alpha(u_0) a'(u_0)$ for some $\alpha$.

		Consider the space-like  plane passing through  $\text{Re}(s^{\prime}(u_0) )$ and $a^{\prime}(u_0)$. 
		
				Let $  \text{Re}(s^{\prime}(u_0) ) \times_L  a^{\prime}(u_0)) = \beta(u_0)  n_s(u_0)$. Since  $ \text{Re}(s^{\prime})$ is spacelike,  $\text{Re}(s^{\prime}(u_0 )) \times_E  a^{\prime}(u_0) $ is timelike and thus we have $\text{Re}(s^{\prime}(u_0 )) \times_L  a^{\prime}(u_0) $ is also timelike and we can choose $\beta$ such that
		$\langle n_s, n_s\rangle_{L} = - 1$.

		Now we take $d_0^{n_s \prime} = n_s \times_L a^{\prime}$ and  $  \Delta^{n_s}  = d_0^{n_s \prime} \times_E a^{\prime}$.  
		Hence,  the vectors $a^{\prime},\;d_0^{n_s \prime},\;\Delta^{n_s}$ are linearly independent and we can write 
		$$\text{Re}(s^{\prime}) = \alpha_1 a^{\prime} + \alpha_2 d_0^{n_s \prime} + \alpha_3 \Delta^{n_s}.$$ Now if we perform $\times_{L} a^{\prime}$ on both sides of the previous equation and noting that $ d_0^{n_s \prime} \times_L a^{\prime}  = n_s \langle a^{\prime}, a^{\prime} \rangle_L $ (see appendix \ref{crossproduct} (1)), we get 
		$$\beta n_s=\text{Re}(s^{\prime}) \times_L a^{\prime} = \alpha_2  d_0^{n_s \prime} \times_L a^{\prime}  + \alpha_3 \Delta^{n_s} \times_L a^{\prime} = \alpha_2 n_s \langle a^{\prime}, a^{\prime} \rangle_L + \alpha_3 \Delta^{n_s} \times_L a^{\prime}.  $$
			
			In other words, $(\beta - \alpha_2 \langle a^{\prime}, a^{\prime} \rangle_L) n_s = \alpha_3 \Delta^{n_s} \times_L a^{\prime}$.
			
			One can show that  $\Delta^{n_s} \times_L a^{\prime}$ is not timelike (see appendix \ref{crossproduct} (2)).  Since $n_s$ is timelike,  $\alpha_3=0$.

		This shows $\text{Re}(s^{\prime}) \perp_E \Delta^{n_s}$, at $u_0$. We repeat for all $u_0$ such that $\text{Re}(s'(u_0))\neq \alpha(u_0) a'(u_0)$ for some $\alpha$.

	\end{proof}
	
	Since $x_0\in \widetilde{\mathcal S}_p$, and $\widetilde{\mathcal S}_p$ is an open subset (in sup-norm topology of $C^\omega(\Omega,\mathbb C^3)$). Therefore, there is an open ball  $B_\epsilon(x_0)\subset \widetilde{\mathcal S}_p$ around $x_0$ and so $B_\epsilon(x_0)\subset \mathcal{J}$ (as a subset). Here, $\mathcal{J}$ has $n-$topology, and if $V^{\Omega}$ is an open subset of $\mathcal{J}$ containing $x_0$, then it is not very obvious that there is some $\epsilon>0$ such that $B_\epsilon(x_0)$(an open ball in sup-norm topology of $C^\omega(\Omega, \mathbb C^3)$) is contained in $V^{\Omega}$. The coming paragraphs show that such an $\epsilon$ exists.

	Let $x_0$ and $V^{\Omega}$ as above. Since $x_0 \in V^{\Omega} =V^{\Omega}\cap \mathcal{J}= \cup_{n \in \mathcal N} V_n $, where $V_n=V^{\Omega}\cap J_n$, open sets of $J_n$ in sup-norm subspace topology (induced  from $C(\Omega, \mathbb C^3))$. Recall $x_0\in J_n$ for all $n\in \mathcal N$ and $x_0\in V_n$ for all $n\in \mathcal N$ (from Proposition \ref{invofHN}). Therefore, for each $n$, there is a ``sliced" sup-norm-open ball $B^{sl}_{\epsilon_n}(x_0)$ around $x_0$ in $V_n$.

	Let $\Omega = \{ z :|z| <1 \}$ be the unit disc.  We define exhaustion by compact sets of $\Omega$ as follows: $\Omega = \cup_{0<r<1} \overline{\Omega_r}$ where $\Omega_r = \{ z : |z| < r \}$.
	
	We define some useful constants, $K_1, K_2, K_3$.
	
	Let $K = \|x_0\|_\infty + 1$ and  $K_1 (r) = K \frac{2}{1-r^2}$.  
	If $0<r<1$, $|\alpha| < 1$ and $|a| < r$, then  $K \frac{1 - |\alpha|^2}{1-|a|^2} \leq  K_1(r)$. 
	
	Next, if $w=(w_1,w_2,w_3)\in \mathcal U_p$ be  such that   $ |w_i| \leq K_1(r)$, then there exists a constant $K_2(r)$ such that $|\text{Re}(w') \times_L a^{\prime}|_E \leq K_2(r)$. We can take $K_2(r) = \sqrt{2} K_1(r) K_a$ where $\sum_{i=1}^3 |a_i^{\prime}|^ 2 \leq K_a^2$ on $\overline{\Omega_r}$ (see appendix \ref{expression for K2}) and we also have $|\text{Re}(w^{\prime}) \times_E a^{\prime}|_{L}\neq 0\; \text{on}\; \overline{\Omega_r} $. Thus $ 0<C(r)\leq -|\text{Re}(s^{\prime}) \times_L a^{\prime}|_{L}\neq 0\; \text{on}\; \overline{\Omega_r}$. Take $K_3(r)=\frac{K_2(r)}{C(r)}$.

	We define the following set (which is a subset of ${\mathcal N}$ and is  just going to be an indexing set).
	\begin{definition} Let 
	\begin{equation} S:= \{ n: \Omega\to \mathbb C^3: \;\;  \|n\|_\infty \leq K_3(r)\;\text{on}\;  \overline{\Omega_r},\;0<r<1\;\text{and}\;\langle n, n\rangle_L=-1  \}
	\end{equation}
	\end{definition}
	Here, we remark that since $S\subset C^\omega(\Omega, \mathbb C^3)$, it has a subspace topology induced from  $C^\omega(\Omega, \mathbb C^3)$.On $S$, we will not take the subspace topology.  Instead, we consider the compact open topology on $S$ and will show it is compact. This can be seen as follows.
	\begin{prop}
	The set $S$ is compact in the compact open topology. 
	\end{prop}
	 
	\begin{proof}
	Let $\{n^l\}_{l=1}^{\infty}$ be a sequence in $S$. Then the $i$-th component of $n^l$, i.e. $n_i^l$ is bounded by $K_3(r)$ for all $z \in \overline{\Omega_r}$. Thus $n_1^l$ is a locally bounded family and thus is normal (by Montel's theorem page 153, Conway (\cite{conway}), and hence there is a sub-sequence which converges to a holomorphic function $n_{01}$. By renaming this sub-sequence by $l$ again, $n_2^l$ is a locally bounded family with a sub-sequence that converges to holomorphic $n_{02}$. Again, renaming this sub-sequence by $l$, we get a locally bounded family $n_3^l$, which has a sub-sequence that converges to $n_{03}$. 
	
	Thus the sequence $n^l $ in $S$ has a limit $n_0=(n_{01},n_{02},n_{03})$ in $C^{\omega}(\Omega, \mathbb C^3)$. It is easy to check $n_{01}^2 + n_{02}^2 - n_{03}^2 =-1$ and 
	for  $ z \in \overline{\Omega_r}$,  $|n_{0i}(z)|  \leq K_3(r) $. Thus $n_0 \in S$. 
	\end{proof}
	
	We take $S$ as an indexing set and define ${ B^{sl}} = \cup_{n \in S} B^{sl}_{\epsilon_n}(x_0)$ (the superscript ``sl''  stands for slice). Take  $\epsilon_0 = {\rm inf}_{ S} \epsilon_n $ and the claim is $\epsilon_0\neq 0$. Suppose $\epsilon_0=0$; then we take a minimizing sequence
	$n_i \in S $ such that $\epsilon_{n_i} \rightarrow 0$. Since $S$ is compact, $n_i \rightarrow n_{\tau}$ and $n_{\tau} \in S$. Consider $ V_{n_{\tau}}$, which is non-empty (follows from Proposition \ref{invofHN}).   Let $B^{sl}_{\epsilon_{n_{\tau}}}(x_0) \subset V_{n_{\tau}}$ be the ``slice'' ball about $x_0$. By assumption $0=\epsilon_0 = \epsilon_{n_{\tau}} >0$, we have a contradiction. This shows that $\epsilon_0 \neq 0$.

	We choose $\epsilon_1>0$ such that $\epsilon_1<\epsilon_0$., and take sup-norm ball $B_{\epsilon_1}(x_0)$ in $\widetilde{\mathcal Sp}$.  By Proposition \ref{SptildemathcalJ} we have  for $s  \in B_{\epsilon_1}(x_0) $,  $s \in J_n$ for some $n \in {\mathcal N}$.

	\begin{prop}\label{nbhd}
		Let $x_0 = a+id_0 \in V^{\Omega}\subset \mathcal{J}$ be $L$-isotropic. 
		Then there is a  sup-norm ball containing $x_0$ in $V^{\Omega}$.  
	\end{prop}
	\begin{proof}
		For $V^{\Omega}$, we take $ B_{\epsilon_1}(x_0)$ as above.  We will  show that $ B_{\epsilon_1}(x_0) \subset B^{sl}_{\epsilon_{n_s}}(x_0)$ for some $n_s \in {\mathcal N}$.  
		
		Let $s\in B_{\epsilon_1}(x_0) $. We have by the previous argument that $s \in J_{n_{s}}$ for some $n_{s} \in  {\mathcal N}$.  But suppose we are able to prove that such $n_s\in S$, then it will prove that  $|| s - x_0||_{\infty} < \epsilon_1  < \epsilon_{0}<\epsilon_n $ for all $n \in S$ and this will prove that  $s \in V_n$ for all $n \in S \subset {\mathcal N}$. Thus $s \in V^{\Omega} = \cup_{n \in {\mathcal N}} V_n$.  
		
		We now show that $n_s\in S$. 
		Let $s=(s_1, s_2, s_3)$. Then $|s_i| < \|x_0\|_{\infty} + \epsilon < \|x_0\|_{\infty} + 1 =K$ for $i=1,2,3$. Then $\frac{|s_i|}{K} < 1 $. Let $a\in \overline{\Omega_r}$ and  $\frac{s_i (a)}{K} = \alpha$ where $|\alpha| <1$.  Then by a consequence of Schwarz's Lemma (page 132 (2.3), Conway),  we have $|s_i^{\prime}(a)| < K \frac{1 - |\alpha|^2}{1-|a|^2} \leq K_1(r)$.  We can repeat for every $a \in \overline{\Omega_r}$ and each $i=1,2,3$. It follows that we have for each $i=1,2,3$, $|s_i^{\prime}| \leq K_1(r)$  on $\overline{\Omega_r}$. Recall since $s \in J_{n_s}$, $\text{Re}(s^{\prime}) $ is spacelike and thus  $\text{Re}(s^{\prime} ) \times_E  a^{\prime} $  is timelike and so is $\text{Re}(s^{\prime} ) \times_L  a^{\prime} $ . Thus  there exists a unit timelike vector $n_s$ and $ \beta: \Omega \rightarrow {\mathbb C}$ such that $\beta n_s = \text{Re}(s^{\prime} ) \times_L  a^{\prime} $. This implies $|\beta|=-|\text{Re}(s^{\prime}) \times_L a^{\prime}|_{L}$.
		
		Now by taking $w=s$ in the definition of the constant $K_2(r)$, we get  $ |n_{s}|_E \leq K_3(r)$. This is because  
		$|\beta||n_s|_E=|\text{Re}(s^{\prime}) \times_L a^{\prime}|_{E}$ and this gives $$|n_s|_E=-\frac{|\text{Re}(s^{\prime}) \times_L a^{\prime}|_{E}}{|\text{Re}(s^{\prime}) \times_L a^{\prime}|_{L}}\leq \frac{K_2(r)}{C(r)}\leq K_3(r).$$
		Thus $|n_{si}|\leq|n_s|_{E} \leq K_3(r)$ for $i=1,2,3$ on $\overline{\Omega_r}$ for $0<r<1$ and this imply $\|n\|_{\infty}\leq K_3(r)$.  Since $s \in J_{n_s}$ implies $n_s \in {\mathcal N}$ by definition, i.e. $n_{s1}^2 + n_{s2}^2 - n_{s3}^2 = -1$, we have  $n_s \in S$.

		{\bf Remark:}
		We have defined various topologies on different spaces,  like sup-norm, compact-open, $n$-topology etc.  Each one has it's role to play and we donot need comparisons of these topologies.

	\end{proof}	
	\color{black}
	Now we come to the main theorem of this section.
	
	\subsection{Local interpolation by maximal surface}
	For $l\in \mathcal{S}_p$ such that for all $z$, $|a^\prime\times l^\prime|\neq 0$ we define 
	\begin{equation}\label{dl_formainthm}
		d^{l}(z) =\int_{u_0}^z\left( \frac{(a^{\prime} \times_L l^{\prime})\times_L l^{\prime}}{ | a^{\prime} \times_L l^{\prime}|}\right)(w) dw. 
	\end{equation}
	$d^l\in {\mathcal W}_{(\Omega,\mathbb C^3, u_0)}$. 
	We state the theorem below. 
	\begin{theorem}\label{mainthm}Let $\Omega$ be as in Notations (\ref{notation1}).  Let $a, d_0 \in \mathcal{ S}_p $ be fixed such that $x_0= a+ i d_0$ is $L$-isotropic and  $X_0= \text{Re}(a+id_0)$ is a generalised maximal immersion in ${\mathbb L}^3$.  
		
		Then there are $\eta_1>0$, and $\eta_2 >0$, such that for all curves $l \in {\mathcal S_p}$  such that 
		\begin{enumerate}
			\item 	$\|l- a\|_\infty<\eta_1, \|l^{\prime}- a^{\prime}\|_\infty<\eta_1$ and $
			\|d^l - d_0\|_{\infty}<\eta_2, \|d^{l\prime} - d_0^{\prime}\|_{\infty}<\eta_2,
			$
			\item For all $z\in \Omega$, $a^\prime(z)\times_L l^\prime (z)\neq 0$ .
		\end{enumerate}
		
		Then there exists a  generalised maximal surface $ X: \Omega \rightarrow \mathbb{L}^3 $ and an automorphism of $\Omega$, $\gamma:\Omega \rightarrow \Omega$, 
		an automorphism such that for all  $u\in I\subset \Omega$, 
		$$X(u,0)=l(u)\;\text{and}\; (X\circ\gamma^{-1})(u,0)=a(u).$$

	\end{theorem}
	
	\begin{proof}
		Let $l\in\mathcal {S}_p$. We take $n_l = \frac{(a^{\prime} \times_L l^{\prime})}{ | a^{\prime} \times_L l^{\prime}|}\in {\mathcal N}$. 
		As before, 
		$x_0= (a + i d_0) \in J_{n}$ all $n\in \mathcal{N}$.  	Let $V^{\Omega} \subset {\mathcal J}$ be an open neighbourhood  (in $n$-topology) of $x_0$ as in Corollary \ref{IFT}. Let $\epsilon>0$ be such that  $B_{\epsilon}(x_0)   \subset V^{\Omega}\subset {\mathcal J}$  as in  Proposition \ref{nbhd}. We take $\eta_1, \eta_2$ be such that $\eta_1^2 + \eta_2^2 < \epsilon$
		
		Let $V_1 = \{s \in {\mathcal L_{(I,\mathbb L^3)}} : \|s-a\|_\infty < \eta_1,\;\; \|s^{\prime}-a^{\prime}\|_\infty < \eta_1   \}$, and  $V_2 = \{ s \in {\mathcal L_{(I,\mathbb L^3)}} : \|d^{s} - d_0\| < \eta_2, \|d^{s \prime} - d_0^{\prime}\| < \eta_2  \}$.  	Then 	${\mathcal S_p} \cap V_1 \cap V_2 $ is an open set in ${\mathcal L_{(I,\mathbb L^3)}}$ since ${\mathcal S_p}$, $V_1$ and $V_2$ are open in ${\mathcal L_{(I,\mathbb L^3)}}$.

		We define a map on $\Omega$ as 
		\begin{equation}C (w) = l(w) +  i  d^l(w) =(C_{1},C_{2},C_{3})(w).  \end{equation}  Since $d^{l \prime}(u)  \perp_L  l^{\prime}(u)$ and $\langle d^{l \prime}(u), d^{l \prime}(u)\rangle_L  = \langle l^{\prime}(u),  l^{\prime}(u)\rangle_L $ for $u \in I$ and hence $\langle C^{\prime}, C^{\prime}\rangle_L(u) =0$ on $I$. By analyticity of  $\langle C^{\prime}, C^{\prime}\rangle_L$ it is zero in $\Omega$. Therefore $C(w)$ defines an $L$-isotropic curve on $\Omega$.

		For $u\in I$, let $\text{Re} (C(u))  = l(u)$.  	Now we claim that $C\in B_{\epsilon}(x_0)  \subset V^{\Omega}$,	$C(w) - (a(w) +id_0(w))  = (l(w) -a(w) ) + i (d^{l} (w) - d_0(w))$. Thus 
		\begin{eqnarray*}
			\|C(w) - (a+id_0)(w)\|_{\infty}^2 && < \| l(w) - a(w) \|_{\infty}^2 + \| d^{l}(w) - d_0(w)\|_{\infty}^2 \\&&< \eta_1^2 + \eta_2^2 <  \epsilon  \text{ since } l \in V_1 \cap V_2 
		\end{eqnarray*}

		Thus $C(w)   \in  B_{\epsilon}(p)   \subset  V^{\Omega} \subset {\mathcal J}$ and hence we can apply  the Corollary ~\ref{IFT}   to  show that      $C = (a+id) \circ \gamma $. Since $C$ is  $L$-isotropic, it is easy to see that $(a+id)$ is also $L$-isotropic.

		We have to show that $X(u,v) = \text{Re} ((a+id) \circ \gamma ) = \text{Re}(C(w))$ is a space like immersion, i.e. we would have a    maximal surface which contains  $l(u) = \text{Re}(C(u)) =X(u,0)$. Notice  $(X \circ \gamma^{-1})(u,v) =(a+id)(u,v)$ contains $a(u)$ when we put $v=0$.  Thus $X(u,0) = l(u)$ and $(X\circ \gamma^{-1})(u,0) = a(u)$, hence $X$ interpolates between $l$ and $a$.

	\end{proof}
	
	\section{{\bf Interpolation by minimal surfaces}}\label{IMiS}
	
	In this section we give an analogous result concerning minimal surfaces.  Let $\Omega$ be as before the open unit disk containing $I$ and let
	 $ {\mathcal E} = {\mathcal E}_{I, \mathbb{E}^3}:=\{e\in C^{\omega}(\Omega,\mathbb{C}^3) \mid: e\;\text{and}\;e'\;\text{are real on}\; I \}$ and ${\mathcal N}':=\{n\in C^{\omega}(\Omega,\mathbb C^3):n|_{I} \; \text{is real},\;\langle n, n\rangle_{E}=1,\; \langle n, a'\rangle_{E} =0\}$. 
	\begin{theorem}\label{minthm}
	Let $\Omega$ be as above and
		let  $a, d_0 \in {\mathcal E}$  be fixed such that $a -i  d_0$ is $E$-isotropic in the sense that $\langle (a-id_0)', (a-id_0)' \rangle_{E}=0$ and  $X_0= \text{Re}(a-id_0)$ is a generalised minimal immersion in ${\mathbb E}^3$.  
		
		Then there are $\eta_1>0$, and $\eta_2 >0$, such that for all  real analytic curves $e\in \mathcal{E}$  such that 
		\begin{enumerate}
			\item 	$\|e- a\|_\infty<\eta_1, \|e^{\prime}- a^{\prime}\|_\infty<\eta_1 \;\; \text{and}\;\;
			\|d^e- d_0\|_{\infty}<\eta_2, \|d^{e\prime} - d_0^{\prime}\|_{\infty}<\eta_2,
			$
			\item For all $z\in \Omega$, $a^\prime(z)\times_E e^\prime (z)\neq 0$ 
		\end{enumerate}

	Then there exists a  generalised minimal surface $ X: \Omega \rightarrow \mathbb{E}^3 $ and an automorphism of $\Omega$,  $\gamma:\Omega \rightarrow \Omega$ such that
		$$X(u,0)=e(u)\;\text{and}\; (X\circ\gamma^{-1})(u,0)=a(u)$$
		where $u\in I\subset \Omega$.
		
	\end{theorem}
	
	\begin{proof}
Much of the proof goes very similar to the maximal surface case.  There is a lot of simplification and all norms and cross products  are Euclidean.We give  a rough sketch  of the proof here. 

First point of difference is we use the Bj\"orling problem for minimal surfaces and the Schwarz solution to it.  Consider $n: I\to \mathbb E^3$ as a real analytic map such that $ \langle n(u), a'(u) \rangle_E =0$ for all $u\in I $ and $\langle n(u),n(u)\rangle_E = 1$ for all $u\in I$.  The map $n$ is referred to as the unit normal vector field along $a$.  The Bj\"orling problem poses the question of whether there exists a minimal  surface $X : \Omega \to\mathbb E^3$ with $I\subset \Omega\subset \mathbb{C}$, where $\Omega$ is simply connected,  and the surface $X$ satisfies the following conditions: $X(u, 0) = a(u)$ and $N(u, 0) = n(u)$ for all $u \in I$.  Here, $N: \Omega\to \mathbb E^3$ is a unit normal to the surface $X$. The well known Schwarz solution to Bj\"orling problem is given by 
\begin{equation}
		X(u,v)=\text{Re}\{a(w)-i\int_{u_0}^{w} n(\widetilde{w})\times_E a'(\widetilde{w})d\widetilde{w}\}; w=u+iv, u_0\in I,
	\end{equation}
	where $n(w)$ and $a(w)$ are the complex analytic extensions of $n(u)$ and $a(u)$ over $\Omega$.
  On $\Omega$, we define 
	\begin{equation}
		d_0^n(z)= \int_{u_0}^z  (n \times_E a^{\prime})(w) dw.
	\end{equation}
	We see $d_0^n(u_0) =0$, and its restriction on $I$ is real. Therefore  $d_0^n\in \mathcal{W}_{(\Omega, \mathbb C^3, u_0)}$.  
	
	Moreover, $$\langle d_0^{n \prime}, d_0^{n \prime} \rangle_E = \langle a^{\prime}, a^{\prime}\rangle_E,\;\; \langle d_0^{n \prime}, a^{\prime}\rangle_E =0.$$    
	In other words, $\langle (a-i d_0^n)^{\prime}, (a-id_0^n)^{\prime}\rangle_E =0$.	We define the map $\Delta^n: I\to \mathbb R^3$ as 
	\begin{equation}\label{deltan}
		\Delta^n(u) := d_0^{n \prime}(u) \times_E a^{\prime}(u).
	\end{equation}
	
	\begin{definition} We have the following
		\begin{enumerate}
			\item $J_{n}^{'}:= \{ c \in  C^{\omega}(\Omega, {\mathbb C}^3)| \text{Re} (c^{\prime}(u))  \perp_E \Delta^n(u)  \; \forall \; u \in I \}$.  
			
			\item ${\mathcal J}^{'} := \cup_{n  \in \mathcal N} J^{'}_{n}$ with the topology as follows:  $V\subset C^{\omega}(\Omega, {\mathbb C}^3)$ is open in ${\mathcal J}^{'}$ if and only if for all $n$, $V \cap J^{'}_{n}$ is open in $J_{n}^{'}$ in sup-norm topology. We  call this $n$-topology.  
		\end{enumerate}
	\end{definition}

	The modication to section 2 is to invert the Bj\"{o}rling-Schwarz  function $$ \tilde{H}:C^{\omega}(\Omega,\Omega)\times  C^{\omega}( \Omega, {\mathbb C}^3) \rightarrow   C^{\omega}(\Omega, {\mathbb C}^3)$$ defined by 
		\begin{equation}
			\tilde{H}(\gamma, d)= (a-id)\circ \gamma
		\end{equation}
		
		Our main aim is to invert the map $\tilde{H} $.  For this we compute the differential $D\tilde{H}$ of $\tilde{H}$,  restrict it to $T_n ':= \{ (V, d) \in {\mathcal C}_{(\Omega,\mathbb C, u_0)} \times {\mathcal W_{(\Omega,\mathbb C^3, u_0)}} :  \text{Re} ((a^{\prime}(u) - i d^{n \prime}_0(u)  ) V(u) - i d(u))^{\prime} \perp_E \Delta^n(u),  u\in I\}$.   For $n\in \mathcal N'$, we define ${\tilde{\mathcal H}}_n : T_n' \rightarrow J^{'}_{n}  $ to be the  map given by ${\tilde{\mathcal H}}_n(V, d) = (a^{\prime}-  i d^{n \prime}_0 ) V - i d$. We show that $\tilde {\mathcal H}_n|_{T_n'}=D\tilde{H}|_{t_n'}$ in the same way as in the maximal surface case. We show that $\tilde{\mathcal H}_n: T_n  \rightarrow J^{`}_{n}$ is invertible.  Then invoking the inverse function theorem for Banach spaces we show that the 
Bj\"{o}rling-Schwartz  function $H$ is invertible in suitable domains. This is completely analogous to the maximal surface case.

 Next we show the existence of sup-norm open ball containing an appropriate $x_0$ in $V^{\Omega}\subset \mathcal{J}^{'}$. 
 Let $x_0 =a - id_0$ be a special point such that $a , d_0 \in {\mathcal E}$ and $\langle a^{\prime}, a^{\prime}\rangle_E = \langle d^{\prime}_0,d^{\prime}_0\rangle_E,  \langle a^{\prime},d^{\prime}_0\rangle_E = 0$.
 Before proving the existsence of  a sup-norm ball containing $x_0$, we prove the following proposition
	\begin{prop}\label{E}
	$ {\mathcal E} \subset \mathcal{J}^{'}$.
	\end{prop}
	
	\begin{proof}
		
		Let $s \in {\mathcal E} $. We claim that $s \in J^{'}_{n_s}$ for some $n_s \in {\mathcal N'}$. For this we just need to show for $u\in I,\;\text{Re}(s^{\prime}(u) ) \perp_E \Delta^{n_s}(u)$. 
		
		If $\text{Re}(s'(u))=\alpha(u) a'(u)$ for all $u \in I $, then $s\in J^{'}_n$ for all $n$, and we are done.
		
		Now, if the above is not true, then there will be some $u_0 \in I$, $\text{Re}(s'(u_0))\neq \alpha(u_0) a'(u_0)$ for some $\alpha$.   
		Consider the   plane passing through  $\text{Re}(s^{\prime}(u_0) )$ and $a^{\prime}(u_0)$. 
	Let $  \text{Re}(s^{\prime}(u_0)  \times_E  a^{\prime}(u_0)) = \beta(u_0)  n_s(u_0)$.  We can choose $\beta$ such that that
		$\langle n_s, n_s\rangle_{E} =  1$. 
		Now we take $d_0^{n_s \prime} = n_s \times_E a^{\prime}$ and  $  \Delta^{n_s}  = d_0^{n_s \prime} \times_E a^{\prime}$.  
		Hence,  the vectors $a^{\prime},\;d_0^{n_s \prime},\;\Delta^{n_s}$ are linearly independent and we can write 
		$$\text{Re}(s^{\prime}) = \alpha_1 a^{\prime} + \alpha_2 d_0^{n_s \prime} + \alpha_3 \Delta^{n_s}.$$ Now if we perform $\times_{E} a^{\prime}$ on both sides of the previous equation and noting that $ d_0^{n_s \prime} \times_E a^{\prime}  = n_s \langle a^{\prime}, a^{\prime} \rangle_E $  we get 
		$$\beta n_s=\text{Re}(s^{\prime}) \times_E a^{\prime} = \alpha_2  d_0^{n_s \prime} \times_E a^{\prime}  + \alpha_3 \Delta^{n_s} \times_E a^{\prime} = \alpha_2 n_s \langle a^{\prime}, a^{\prime} \rangle_E + \alpha_3 \Delta^{n_s} \times_E a^{\prime}.  $$
			
			In other words, $(\beta - \alpha_2 \langle a^{\prime}, a^{\prime} \rangle_E) n_s = \alpha_3 \Delta^{n_s} \times_E a^{\prime}$. Since RHS is proportional to $d_0^{n_s \prime} $ and $n_s$ is perpendicular to 
			$d_0^{n_s \prime} $, $\alpha_3=0$.

		This shows $\text{Re}(s^{\prime}) \perp_E \Delta^{n_s}$ at $u_0$.  We repeat for all $u_0$ such that $\text{Re}(s'(u_0))\neq \alpha(u_0) a'(u_0)$ for some $\alpha$.

	\end{proof}

	Next propostion is
	
	\begin{prop}
	 ${\mathcal N'}$ is compact in the compact open topology.
	\end{prop}
	
	\begin{proof}
	This is easy to show since $n \in {\mathcal N'}$ has to have $n_1^2 + n_2^2 + n_3^3 =1$ where $n =(n_1, n_2, n_3)$. 
	\end{proof}
		Since $x_0=a-id_0$ is $E$-isotropic, $x_0\in \mathcal E $ and ${\mathcal E}$ is an open subset (in sup-norm topology of $C^\omega(\Omega,\mathbb{ C}^3)$). Therefore, there is an open ball  $B_\epsilon(x_0)\subset \mathcal{E}$ around $x_0$ and so $B_\epsilon(x_0)\subset  \mathcal{J}^{'}$. Here,  $ \mathcal{J}^{'}$ has $n-$topology, and if $V^{\Omega}$ is an open subset of  $\mathcal{J}^{'}$  containing $x_0$, then it is not very obvious that there is some $\epsilon>0$ such that $B_\epsilon(x_0)$(an open ball in sup-norm topology $C^\omega(\Omega,\mathbb{ C}^3)$) is contained in $V^{\Omega}$.  We show that such an $\epsilon$ exists. 
	We take ${\mathcal N'}$ as an indexing set and define ${ B^{sl}} = \cup_{n \in {\mathcal N'}} B^{sl}_{\epsilon_n}(x_0)$ (the superscript ``sl''  stands for slice). Take  $\epsilon_0 = {\rm inf}_{ {\mathcal N'}} \epsilon_n $ and the claim is $\epsilon_0\neq 0$. Suppose $\epsilon_0=0$; then we take a minimizing sequence
	$n_i \in S $ such that $\epsilon_{n_i} \rightarrow 0$. Since ${\mathcal N}$ is compact, $n_i \rightarrow n_{\tau}$ and $n_{\tau} \in \mathcal N'$. Consider $ V_{n_{\tau}}$, which is non-empty.    Let $B^{sl}_{\epsilon_{n_{\tau}}}(x_0) \subset V_{n_{\tau}}$ be the ``slice'' ball about $x_0$. By assumption $0=\epsilon_0 = \epsilon_{n_{\tau}} >0$, we have a contradiction. This shows that $\epsilon_0 \neq 0$.

	We choose $\epsilon_1>0$ such that $\epsilon_1<\epsilon_0$., and take sup-norm ball $B_{\epsilon_1}(x_0)$ in $\mathcal E$.  By Proposition \ref{E} we have  for $s  \in B_{\epsilon_1}(x_0) $,  $s \in J^{'}_n$ for some $n \in {\mathcal N'}$. 
	
	\begin{prop}
		Let $x_0 = a-id_0 \in V^{\Omega}\subset  \mathcal{J}^{'}$ be $E$-isotropic. Then there is a  sup-norm ball containing $x_0$ in $V^{\Omega}$.  
	\end{prop}
	
		Rest of the proof of  Theorem \ref{minthm} goes as follows:
	
	Let $e$ as in the theorem. We take $n_e = \frac{(a^{\prime} \times_E e^{\prime})}{ | a^{\prime} \times_E e^{\prime}|}\in {\mathcal N'}$. 
		Now, observe that
		$x_0= (a - i d_0) \in J_{n}^{'}$ all $n\in \mathcal{N'}$.  	Let $V^{\Omega} \subset {\mathcal J}^{'}$ be an open neighbourhood  (in $n$-topology) of $x_0$. Let $\epsilon>0$ be such that  $B_{\epsilon}(x_0)   \subset V^{\Omega}\subset \mathcal{J}^{'}$ . We take $\eta_1, \eta_2$ be such that $\eta_1^2 + \eta_2^2 < \epsilon$
		
		Let $V_1 = \{s \in {\mathcal E} : \|s-a\|_\infty < \eta_1,\;\; \|s^{\prime}-a^{\prime}\|_\infty < \eta_1   \}$, and  $V_2 = \{ s \in {\mathcal E} : \|d^{s} - d_0\| < \eta_2, \|d^{s \prime} - d_0^{\prime}\| < \eta_2  \}$.  Since ${\mathcal E}$, $V_1$ and $V_2$ are open in ${\mathcal E}$, the set	${\mathcal E} \cap V_1 \cap V_2 $ is an open set in ${\mathcal E}$

		We define a map on $\Omega$ as 
		\begin{equation}C (w) = e(w) -  i  d^e(w) =(C_{1},C_{2},C_{3})(w).  \end{equation}  Since $d^{e \prime}(u)  \perp_E  e^{\prime}(u)$ and $\langle d^{e \prime}(u), d^{e \prime}(u)\rangle_E  = \langle e^{\prime}(u),  e^{\prime}(u)\rangle_E $ for $u \in I$ and hence $\langle C^{\prime}, C^{\prime}\rangle_E(u) =0$ on $I$. By analyticity of  $\langle C^{\prime}, C^{\prime}\rangle_E$ it is zero in $\Omega$. Therefore, $C(w)$ defines an $E$-isotropic curve on $\Omega$.

		For $u\in I$, let $\text{Re} (C(u)) = e(u)$.  	Now we claim that $C\in B_{\epsilon}(x_0)  \subset V^{\Omega}$,	$C(w) - (a(w) -id_0(w))  = (e(w) -a(w) ) + i (d^{e} (w) - d_0(w))$. Then 
		\begin{eqnarray*}
			\|C(w) - (a+id_0)(w)\|_{\infty}^2 && < \| e(w) - a(w) \|_{\infty}^2 + \| d^{e}(w) - d_0(w)\|_{\infty}^2 \\&&< \eta_1^2 + \eta_2^2 <  \epsilon  \text{ since } e\in V_1 \cap V_2 
		\end{eqnarray*}

		Thus $C(w)   \in  B_{\epsilon}(p)   \subset  V^{\Omega} \subset {\mathcal J}^{'}$ and hence see to  that      $C = (a-id) \circ \gamma $. Since $C$ is  $E$-isotropic, it is easy to see that $(a-id)$ is also $E$-isotropic.

		We have to show that $X(u,v) = \text{Re} ((a-id) \circ \gamma ) = \text{Re}(C(w))$ is an immersion, i.e. we would have a   minimal  surface which contains  $e(u) = \text{Re}(C(u)) =X(u,0)$. Notice  $(X \circ \gamma^{-1})(u,v) =(a-id)(u,v)$ contains $a(u)$ when we put $v=0$.  Thus $X(u,0) = e(u)$ and $(X\circ \gamma^{-1})(u,0) = a(u)$, hence $X$ interpolates between $e$ and $a$.

	\end{proof}

	\section{{\bf Future Directions}} 
	
	1.  Our methods can be extended to interpolation of two real analytic curves by Born Infeld solitons  and also by timelike minimal surfaces 
	
	2.  We can address the question of constructing maximal surfaces with singularities using interpolation methods.
	
	\section{\bf{Appendix}}
	
	\subsection{$\tilde{d}(u,0)$ is real}\label{d(u,0)realonI}  In the Proposition \ref{invertible} we claimed that  $\tilde{d}(u,0)$ is real. Below we prove it. 
	When $\Delta_3^n(u)\neq 0$, we get 
	
	\begin{equation}
		\begin{split}
			\text{Re} (\tau_3(u)) =  \text{Re} (C(u) \widetilde{V}(u)) = (a_3' \widetilde{V}_1   - d_{03}^{n \prime} \widetilde{V}_2)(u) \\
			\text{Re} (\tau_3(u) ) = a_3' \frac{ \text{Re} \tau_2 d_{01}^{n \prime}- \text{Re} \tau_1 d_{02}^{n \prime}}{ \Delta^n_3}(u) + d_{03}^{n \prime} \frac{\text{Re} \tau_1 a_2'- \text{Re} \tau_2 a_1^{\prime}}{\Delta^n_3}(u).
		\end{split}
	\end{equation}
	
	In other words,
	
	\begin{equation}
		\begin{split}
			\Delta^n_3(u) \text{Re} (\tau_3(u)) = (- \Delta^n _1 \text{Re} \tau_1 - \Delta^n_2 \text{Re} \tau_2)(u) 
		\end{split}
	\end{equation}
	where on $I$, $\Delta^n_1 =  a_3'd_{02}^{n \prime} - a_2'd_{03}^{n \prime}$, 
	$\Delta^n_2 =  a_1'd_{03}^{n \prime}- a_3'd_{01}^{n \prime} $ and  $\Delta^n_3 = a_2'd_{01}^{n \prime} -a_1'd_{02}^{n \prime}$ are the components of $\Delta ^n= d_0^{n \prime} \times_E a^{\prime}$.
	
	In other words, $(\Delta^n_3 \text{Re}(\tau_3) + \Delta^n_2 \text{Re}(\tau_2)  + \Delta^n_1 \text{Re} (\tau_1))(u) =0$ for this $u \in I$. 
	
	Note that this equation is symmetric in all coordinates and is satisfied by our choice  of $\tau = s^{\prime}$.

	Thus ${\widetilde{d}}(u,0)$ is real on $I$.  It is unique as is clear from the equations.

	\subsection{Definition of $V_1$ is consistent}\label{Vconsistent}
	
	In Proposition \ref{invertible} we defined  on $I$
	$$\widetilde{V}_1(u)=\frac{ \text{Re} \tau_2 (u) d_{01}^{n\prime}(u)- \text{Re} \tau_1(u) d_{02}^{n \prime}(u)}{\Delta^n_3(u)},$$ when $\Delta^n_3(u) \neq 0$.
	
	If $\Delta_1^n(u) \neq 0$, it can be defined as 
	
	$$\widetilde{V}_1(u)=\frac{ \text{Re} \tau_3 (u) d_{02}^{n\prime}(u)- \text{Re} \tau_2(u) d_{03}^{n \prime}(u)}{\Delta^n_1(u)},$$
	and if $\Delta_2^n(u) \neq 0$ it can be defined as 
	
	$$\widetilde{V}_1(u)=\frac{ \text{Re} \tau_1 (u) d_{03}^{n\prime}(u)- \text{Re} \tau_3(u) d_{01}^{n \prime}(u)}{\Delta^n_2(u)},$$ when $\Delta^n_2(u) \neq 0$.
	
	Similar definitions for $\widetilde{V}_2$ can be given.
	
	To show that the first two definitions of $\widetilde{V}_1$ are consistent, we observe the following.
	
	Since $\tau = s'$ and $\text{Re}(\tau) \perp_E \Delta^n$ on $I$, we have $\text{Re}\tau_2 \Delta^n_2 + \text{Re} \tau_3 \Delta^n_3 + \text{Re}\tau_1 \Delta_1^n =0$. 
	
	Multiplying the above expression by $d_{02}^{n\prime}$ we get 
	$$\text{Re} \tau_2 \times (-d_{02}^{n\prime} \Delta^n_2) = d_{02}^{n\prime} ( \text{Re}\tau_1 \Delta^n_1 + \text{Re} \tau_3 \Delta^n_3 ).$$  Using the fact that $d_{0}^{n \prime} \perp_E \Delta^n$, and above equality,  we have $$\text{Re} \tau_2 ( d_{01}^{n \prime} \Delta^n_1 + d_{03}^{n \prime} \Delta_3^n)=\text{Re} \tau_2 (-d_{02}^{n\prime} \Delta^n_2) = d_{02}^{n\prime} ( \text{Re}\tau_1 \Delta^n_1 + \text{Re} \tau_3 \Delta^n_3 ).$$
	
	Or, in other words,
	$$\text{Re}\tau_2 d_{01}^{n \prime} \Delta_1^n -  \text{Re}\tau_1 d_{02}^{n \prime} \Delta_1^n = \text{Re}\tau_3 d_{02}^{n \prime} \Delta_3^n - \text{Re}\tau_2 d_{03}^{n \prime} \Delta_3^n.$$
	This proves-
	$$\frac{\text{Re}\tau_2 d_{01}^{n \prime}  -  \text{Re}\tau_1 d_{02}^{n \prime}}{\Delta_3^n} = \frac{\text{Re}\tau_3 d_{02}^{n \prime}  - \text{Re}\tau_2 d_{03}^{n \prime}}{ \Delta_1^n}.$$

	\subsection{Definition of $\widetilde{V}$ on all of $\Omega$:}\label{defn of tilde V on all of Omega}
	First, we  show that on $\Omega$, $\Delta^n \neq 0$.
	
	$\Delta^n = d_0^{n \prime} \times_E a^{\prime} = (n \times_L a^{\prime} )\times_E a^{\prime}$.
	Suppose $\Delta^n(z_0) =0$ for $z_0 \in \Delta^n$.
	
	This implies $(n \times_L a^{\prime} \times_E a^{\prime}) (z_0)=0$ and hence $(n \times_L a^{\prime}) \times_L a^{\prime} =0$. 
	
	But $(n \times_L a^{\prime}) \times_L a^{\prime} = < a^{\prime}, a^{\prime} >_L n $.(This can be shown using the fact that $<n, a^{\prime}>_L =0$). 
	
	$< a^{\prime}, a^{\prime} >_L(z_0) n(z_0) =0$. Since $n$ is time-like on $I$, $(n_1^2 + n_2^2 - n_3^2 )(z_0)= -1$. Hence $n(z_0) \neq 0$ and 
	we must have $<a^{\prime}, a^{\prime}>_L (z_0)=0$.  We assume that $\Omega$ is small enough that $<a^{\prime}, a^{\prime}>_L (z) \neq 0$ for all $z \in \Omega$.

	Let us consider the formula 
	$$\widetilde{V}_1(u)=\frac{ \text{Re} \tau_2 (u) d_{01}^{n\prime}(u)- \text{Re} \tau_1(u) d_{02}^{n \prime}(u)}{\Delta^n_3(u)},$$ when $\Delta^n_3(u) \neq 0$.

	Let $T_j(u) = \text{Re}(\tau_j) (u,0)$ where $j=1,2,3$.
	
	$T_j$ is a real analytic function of $u$ and hence has an analytic extension

	If $\Delta_3^n =0 $ for some $z=z_0$ then we can use the formula $$\widetilde{V}_1(z)=\frac{ T_2 (z) d_{01}^{n\prime}(z)- T_1(z) d_{02}^{n \prime}(z)}{\Delta^n_3(z)},$$ when $\Delta^n_3(z) \neq 0$.
	
	If $\Delta^n_3(z_0) =0$ for some $z_0 \in \Omega$, we can use an appropriate  equivalent formula.
	
	The same type of extension holds for $\widetilde{V}_2$.

	\subsection{Solving the ODE}\label{solving ODE}
	
	First  We claim $a^{\prime} + i d_0^{n \prime} \neq 0 $ on $\Omega$. Since $<a^{\prime}, a^{\prime}>_L \neq 0$ on $\Omega$, 
	
	Suppose $a^{\prime} + i d_0^{n \prime} =0$. Then
	$< a^{\prime}, a^{\prime}>_ L + i < a^{\prime}, d_0^{n \prime}>_L =0$. Since $d_0^{n \prime} = n \times_L a^{\prime} $, $n $ being a unit timelike vector, we have 
	$< a^{\prime}, d_0^{n \prime}>_L =0$. Since $<a^{\prime}, a^{\prime}>_L \neq 0$ on $\Omega$, we have a contradiction.
	
	Now we show  the following ODE has a solution. 
	
	$(a^{\prime} + i d_0^{n \prime}) \widetilde{V} = ((a^{\prime} + i d_0^{n \prime}) V)^{\prime}$, $V(u_0)=0$.

	By the method of introducing integrating factors, we have
	on $I$ $$V(u) = \frac{1}{(a^{\prime} + i d_0^{n \prime} )(u)} \int_{u_0}^u (a^{\prime} + i d_0^{n \prime}) \widetilde{V} dt .$$
	This is well defined and can be extended to $\Omega$ since $a^{\prime} + i d_0^{n \prime} \neq 0$ on $\Omega$. 
	
	Details are here- 
	on $I$, 
	$\tilde{V}= \frac{(a^{\prime} + i d_0^{n \prime})^\prime}{ (a^{\prime} + i d_0^{n \prime})}V+ V^\prime$
	
	For some branch of log, we have 
	
	$$\tilde{V}= [\rm{log}(a^{\prime} + i d_0^{n \prime})]^\prime V+ V^\prime$$
	
	We call $F= [\rm{log}(a^{\prime} + i d_0^{n \prime})]^\prime$, and ODE is 
	$$V^\prime+ F V= \tilde{V}$$
	
	We have Integrating factor $$IF=e^{\int F(u) d(u)}=  e^{\rm{log}(a^{\prime} + i d_0^{n \prime})}= (a^{\prime} + i d_0^{n \prime})$$
	
	and we get 
	$$V(u)= \frac{1}{IF}\int IF. \tilde{V} dt$$

	\subsection{Crossproduct}\label{crossproduct}
	
	(1) 
	One can show that $d_0^{n_s \prime} \times_L a^{\prime}  = n_s \langle a^{\prime}, a^{\prime} \rangle_L $ from  a straightforward calculation using the fact that $\langle n_s,  a^{\prime} \rangle =0$.

	
	
	
	
	(2) On $I$,  let $U = \Delta^{n_s} \times_L a^{\prime}=(a,b,c)$ and $V= \Delta^{n_s} \times_E a^{\prime}=(a,b,-c)$.  It can be easily seen that $U$ and $V$ are of the same causal nature.  Since $ \Delta^{n_s}$ is  timelike (being the $E$-cross product of two space like vectors) and $a^{\prime}$ is spacelike,  $U$ can not be timelike, if it was not the case, let $(x,y,z)=\Delta^{n_s}$. Then by reverse Cauchy-Schwarz's inequality (see \cite{Lo}), for these two timelike vectors $U$ and $\Delta^{n_s}$ we have $|\langle U, \Delta^{n_s}\rangle_{L}|\geq \|U\|_{L}\|\Delta^{n_s}\|_L>0 $. This imply $|\langle (a,b,-c),(x,y,z)\rangle_{L}|>0$.  Also, note that $$0<\langle (a,b,-c),(x,y,z)\rangle_{L}=\langle (a,b,c),(x,y,z)\rangle_{E}.$$
	This imply $\Delta^{n_s}$ is not Euclidean perpendicular to $V$, a contradiction.

	\subsection{Expression for $K_2$}\label{expression for K2} In section (3.1) we derived a possible expression for $K_2(r)$ such that 
	if $w=(w_1,w_2,w_3)\in C^\omega(\Omega, \mathbb C^3)$ satisfies    $ \|w_i\|_\infty \leq K_1(r)$, then there exists a constant $K_2(r)$ such that $|\text{Re}(w) \times_L a^{\prime}|_E \leq K_2(r)$. We show that we can  take $K_2(r) = \sqrt{2} K_1(r) K_a$ where $\sum_{i=1}^3 |a_i^{\prime}|^ 2 \leq K_a^2 $ on $\Omega$.
	
	$\text Re(w) \times_L a^{\prime} =(A_1, A_2, A_3) $.
	$A_1 = (\text{w}_2 a_3^{\prime} - w_3 a_2^{\prime}) $ and hence $|A_1|^2 \leq K_1(r)^2 | a_3^{\prime} + a_2^{\prime} |^2 \leq K_1(r)^2 (|a_2^{\prime}|^2 + |a_3^{\prime}|^2 )$.
	
	Similarly, one can work out for $A_2$ and $A_3$.
	
	$|\text Re(w) \times_L a^{\prime}|_E^2 = |A_1|^2 + |A_2|^2 + |A_3|^2 \leq 2 K_1(r)^2 K_a^2 $. Thus we can take $K_2(r) = \sqrt{2} K_1(r) K_a$.

	\section{\bf{Acknowledgement}}
	
	The authors would like  acknowledge the helpful discussions many years ago with Professor Joseph Oesterle and Professor Michael Wolf which helped crystallize the initial formulation of the problem (for minimal surfaces). The authors would also like to thank Dr. Pradip Kumar for his enormous help in reformulating the problem,   calculations and scrutiny.

	This research was supported in part by the International Centre for Theoretical Sciences (ICTS)- Bangalore, for participating in the program - Geometry and Topology for Lecturers (Code: ICTS/gtl2018/06).
	Rukmini Dey would like to acknowledge the support of the  the Department of Atomic Energy,  Government of India under project no. RTI4001.
	Rahul Kumar Singh would like to acknowledge the external grant he has obtained, namely MATRICS (File No. MTR/2023/000990), which has been sanctioned by the SERB.

\end{document}